\documentclass{amsart}

\usepackage[utf8]{inputenc}	

\usepackage{arydshln}

\usepackage[width=.85\textwidth]{caption}

\usepackage{url}

\usepackage[all]{xy}
\usepackage{pb-diagram, pb-xy}

\usepackage{amssymb}
\usepackage{amsthm}

\newcommand{\bigrowheight}{\rule[-0.5em]{0em}{1.7em}}
\newcommand{\medrowheight}{\rule[-0em]{0em}{1em}}

\newtheorem{thm}{Theorem}[section]
\newtheorem{lem}[thm]{Lemma}
\newtheorem{cor}[thm]{Corollary}
\newtheorem{prop}[thm]{Proposition}
\newtheorem{rem}[thm]{Remark}

\newcommand{\bbC}{{\mathbb C}}
\newcommand{\bbF}{{\mathbb F}}
\newcommand{\bbG}{{\mathbb G}}

\newcommand{\bbN}{{\mathbb N}}

\newcommand{\bbZ}{{\mathbb Z}}

\newcommand{\frakg}{{\mathfrak g}}

\newcommand{\fraksl}{{\mathfrak{sl}}}
\newcommand{\frakgl}{{\mathfrak{gl}}}

\newcommand{\frakS}{{\mathfrak S}}
\newcommand{\osp}{{\mathfrak{osp}}}

\newcommand{\Sp}{\mathit{Sp}}

\newcommand{\GO}{\mathit{GO}}
\newcommand{\GSO}{\mathit{GSO}}
\newcommand{\Sl}{\mathit{Sl}}
\newcommand{\Gl}{\mathit{Gl}}

\newcommand{\OSp}{\mathit{OSp}}
\newcommand{\Rep}{{\mathit{Rep}}}

\DeclareMathOperator{\Aut}{Aut}
\DeclareMathOperator{\Out}{Out}

\DeclareMathOperator{\Ind}{Ind}

\newcommand{\Dbc}{{{D^{\hspace{0.01em}b}_{\hspace{-0.13em}c} \hspace{-0.05em} }}}

\newcommand{\End}{{\mathit{End}}}

\newcommand{\one}{{\mathbf{1}}}

\newcommand{\bfD}{{\mathbf{D}}}

\begin{document}

\title[Tannakian categories with a small tensor generator]{Semisimple super Tannakian categories\\ with a small tensor generator}
\author{Thomas Kr\"amer and Rainer Weissauer}
\address{Mathematisches Institut\\ Ruprecht-Karls-Universit\"at Heidelberg\\ Im Neuenheimer Feld 288, D-69120 Heidelberg, Germany}
\email{tkraemer@mathi.uni-heidelberg.de} \email{rweissauer@mathi.uni-heidelberg.de}

\begin{abstract}
We consider semisimple super Tannakian categories generated by an object whose symmetric or alternating tensor square is simple up to trivial summands. Using representation theory, we provide a criterion to identify the corresponding Tannaka super groups that applies in many situations.  As an  example we discuss the tensor category generated by the convolution powers of an algebraic curve inside its Jacobian variety.
\end{abstract}

\maketitle

\thispagestyle{empty}

\section{introduction}

The goal of this paper is to study reductive super groups with a representation whose symmetric or alternating square is irreducible or splits into an irreducible plus a trivial representation. This problem is motivated by the analogous question for general symmetric tensor categories, which are ubiquitous and arise naturally in many areas of mathematics such as representation theory, topology and algebraic geometry \cite{Deligne_HodgeCyclesMotivesShimura}~\cite{GL_Faisceaux}~\cite{KW_Vanishing}~\cite{Kraemer_Semiabelian}. Over an algebraically closed field $k$ of characteristic zero a result by Deligne says that every $k$-linear symmetric tensor category with certain finiteness properties is equivalent to the category $\Rep_k(G, \varepsilon)$ of finite-dimensional algebraic super representations $V$ of an affine super group scheme~$G$ over~$k$, subject to the condition that a certain element $\varepsilon \in G(k)$ acts via the parity automorphism on $V$. This in particular includes the categories of ordinary group representations by taking $\varepsilon = 1$. Since $\Rep_k(G,\varepsilon)$ is a full subcategory of the category $\Rep_k(G)$ of all finite-dimensional algebraic super representations of $G$, it suffices in what follows to consider the categories $\Rep_k(G)$.

\medskip

If such a tensor category has a tensor generator~$X$ in the sense that any object is a subquotient of $(X \oplus X^\vee)^{\otimes r}$ for some $r\in \bbN$, the super group scheme~$G$ is of finite type over~$k$ and is called the Tannaka super group of the category. We then have a faithful algebraic super representation
\[
 G \; \hookrightarrow \; \Gl(V)
\]
on a finite-dimensional super vector space $V$ corresponding to $X$. We say that~$G$ is~{\em reductive} if the category $\Rep_k(G)$ is semisimple. The reductive super groups have been classified in~\cite{Weissauer_SS} and are built up from ordinary reductive groups and the orthosymplectic super groups~$\OSp_{1 | 2m}(k)$. Note that the semisimple case is relevant also for more general tensor categories since a universal construction of Andr\'e and Kahn~\cite[sect.~8]{AK_Nilpotence}~\cite{AK_Erratum} associates to a non-semisimple tensor category a maximal semisimple quotient.  However, the semisimple tensor categories that arise in this way are often hard to approach explicitly. For examples from algebraic geometry we refer to~\cite{KW_Schottky} \cite{Weissauer_BN} \cite{Weissauer_Torelli}; similar examples also arise in representation theory.  
The goal of this paper is to give criteria that allow to identify the Tannaka super group of such semisimple tensor categories if the symmetric or alternating square of a tensor generator splits into only few irreducible pieces. As a typical application we discuss in section~\ref{sec:BN} the tensor category generated by the convolution powers of an algebraic curve in its Jacobian variety; here our criterion considerably simplifies the proof of the main result of~\cite{Weissauer_BN}.

\medskip

Returning to representation theory, let $G$ be a reductive super group over an algebraically closed field $k$ of characteristic zero. For simplicity, in what follows the term {\em representation} refers to a representation on a super vector space if we are dealing with true super groups but to an ordinary representation otherwise.  For~$V$ in~$\Rep_k(G)$ we denote by
\[
 T_\epsilon(V) \;=\; 
 \begin{cases}
  S^2(V) & \textnormal{for $\epsilon = +1$}, \\
  \Lambda^2(V) & \textnormal{for $\epsilon = -1$},
 \end{cases}
\]
the symmetric and the alternating square with respect to the super commutativity constraint of \cite{Deligne_Tensorielles}. If~$T_\epsilon(V)$ is irreducible or a direct sum of an irreducible and a one-dimensional trivial representation $\one$,  we say $V$ is $\epsilon$-{\em small} (or just {\em small}). Small representations are irreducible. If the trivial direct summand $\one$ occurs in $T_\epsilon(V)$, then $V$ is isomorphic to its dual $V^\vee$ and hence carries a non-degenerate symmetric or alternating bilinear form. We say that $V$ is {\em very small} if both~$S^2(V)$ and $\Lambda^2(V)$ are irreducible. Since $\dim_k (\End_G(V\otimes V)) = \dim_k (\End_G(V \otimes V^\vee))$, this is the case iff~$V\otimes V^\vee \cong W \oplus \one$ for some irreducible representation $W\in \Rep_k(G)$.

\medskip

By definition a super group is {\em quasisimple} if it is a perfect central extension of a~(non-abelian) simple super group. For the finite quasisimple groups $G$ very small and self-dual small faithful representations have been classified by Magaard, Malle and Tiep~\cite[th.~7.14]{MMT_Irreducibility}, using earlier results of Magaard and Malle~\cite{MM_Irreducibility}. In a more general setup the list of very small representations has been extended by Guralnick and Tiep to arbitrary reductive groups~\cite[th.~1.5]{GT_Decompositions}. In particular, except for the standard representation of the special linear group, very small representations of~$G$ only exist if the quotient $G/Z(G)$ by the centre $Z(G)$ is finite. The class of small representations is much richer and contains several cases with $\dim(G/Z(G))>0$.

\medskip
To state our main result we use the following notations.
For super groups~$G_i$ and representations $V_i \in \Rep_k(G_i)$ define $G_1\otimes G_2 \subset \Gl(V_1\boxtimes V_2)$ to be the image of the exterior tensor product representation. If a group of automorphisms of~$G_1\otimes G_2$ contains elements that interchange the two subgroups $G_1 \otimes \{ 1\}$ and $\{ 1\} \otimes G_2$, we say that it {\em flips the two factors}. 
If a group acts transitively on a set $X$ and if the action on the set of $2$-element subsets of $X$ is still transitive, we say that the group acts {\em $2$-homogenously} on $X$. If for $V\in \Rep_k(G)$ the restriction $V|_K$ to some normal abelian subgroup~$K\unlhd G$  splits into a direct sum of pairwise distinct characters that are permuted $2$-homogenously and faithfully by the adjoint action of~$G/K$, we say that the representation $V$ is {\em $2$-homogenous monomial}. Finally, a finite $p$-group~$E$ is called  {\em extraspecial} if $E/Z(E)$ is elementary abelian and $Z(E) = [E, E]$ is cyclic of order $p$. Then $|E|=p^{1+2n}$ for some $n\in \bbN$, and for any non-trivial character~$\omega$ of~$Z(E)\cong \bbZ/p\bbZ$ there is a unique irreducible representation $V_\omega\in \Rep_k(E)$ with dimension $p^n$ on which $Z(E)$ acts via $\omega$~\cite[th.~31.5]{Dornhoff}.

\begin{thm} \label{thm:main}
Let $G$ be a reductive super group and $V\in \Rep_k(G)$ an \mbox{$\epsilon$-small} faithful representation of super dimension $d>0$. Then 
one of the following holds: \medskip

\begin{enumerate}
\item The connected component $G^0\subseteq G$ is quasisimple and the restriction $V|_{G^0}$ remains $\epsilon$-small. In this case the possible Dynkin types of $G^0$ and the highest weights of~$V|_{G^0}$ are in table~\ref{tab:smallreps}, using the notations of section~\ref{sec:Lie}. 
\medskip

\item $(G^0, V|_{G^0}) \cong (G_1\otimes G_1, W\boxtimes W)$ where $G_1 \in \{ \Sl_m(k), \Gl_m(k)\}$ and where~$W$ is the $m$-dimensional standard representation or its dual. Here~$G$ flips the two factors so that $G\cong G^0 \rtimes \bbZ/2\bbZ$ and $\epsilon = -1$. \medskip
\item There exists an embedding $G\hookrightarrow \GO_4(k)$ such that $V$ is the restriction of the four-dimensional orthogonal standard representation, and $\epsilon = +1$. \medskip

\item The representation $V$ is $2$-homogenous monomial; then $\epsilon = -1$ unless $V$ has (non-super) dimension $\dim_k(V)\leq 2$. \medskip

\item The group $G = Z(G)\cdot S$ is a (not necessarily direct, but commuting) product of its centre and some finite subgroup $S\subseteq G$. Furthermore we have an exact sequence $0\to H \to S \to \Out(H)$ where
 \smallskip

\begin{enumerate}
\item[(e$_1$)] 
either $H$ is quasisimple, \smallskip

\item[(e$_2$)] 
or $(H, V|_H)\cong (G_1 \otimes G_1, W\boxtimes W)$ for some very small $W\in \Rep_k(G_1)$, in which case $S$ flips the two factors and $\epsilon = -1$, \smallskip

\item[(e$_3$)] or $H$ is a finite $p$-group for some prime $p$ and contains a $G$-stable extraspecial subgroup $E$ of order $p^{2n+1}$ for some $n\in \bbN$. Then~$V|_E$ is irreducible with dimension $p^n$.
\end{enumerate}
\end{enumerate}
\end{thm}

\medskip

By definition of the super commutativity constraint~\cite{Deligne_Tensorielles}, the parity flip $W=\Pi V$ satisfies $S^2(W)=\Lambda^2(V)$ and $\Lambda^2(W)=S^2(V)$ but changes the sign of the super dimension. Furthermore, the super dimension of an irreducible representation of a reductive super group is always non-zero~\cite[lemma~15]{Weissauer_SS}. This explains why we assumed $d>0$ in theorem~\ref{thm:main}.

\medskip

Note that for any faithful irreducible $V\in \Rep_k(G)$, Schur's lemma implies that the centre $Z(G)$ acts on $V$ via scalar matrices. So either $Z(G)=\bbG_m$ or~$Z(G)$ is a finite cyclic group. If the restriction $V|_{G^0}$ to the connected component remains irreducible, then Schur's lemma also applies for the centralizer $Z_G(G^0)\subseteq G$. Thus in the situation of case~(a) the group of connected components is easily controlled since
$
 G/(G^0 \cdot Z_G(G^0)) \; \hookrightarrow \; \Out(G^0)
$
must be a subgroup of outer automorphisms that fix the isomorphism type of the representation $V|_{G^0}$ in table~\ref{tab:smallreps}.

\medskip

For the converse of theorem~\ref{thm:main} one readily checks that all representations $V$ in case (a), (b), (e$_2$) are small. Concerning~(c) recall that the group of orthogonal similitudes $\GO_4(k)$ is the product of its centre with $\GSO_4(k)\cong \Gl_2(k) \otimes \Gl_2(k)$, and that for the latter any small representation must be a product of two very small ones. As a typical example of (d), for any $2$-homogenous subgroup $F$ of the symmetric group $\frakS_d$ we have the \mbox{$2$-homogenous} monomial small representation of~$G = (\bbG_m)^d \rtimes F$ on $V=k^d$ with the natural action. Apart from a single extra case, the $2$-homogenous permutation groups on $d\geq 4$ letters are precisely the doubly transitive ones~\cite[prop.~3.1]{Kantor_AutomorphismGroups} \cite{Kantor_Homogenous}, and the finite doubly transitive groups have been classified by Huppert, Hering and others~\cite[sect.~7.7]{DM_Permutation}. In the extraspecial case (e$_3$) the analysis of the smallness condition is more subtle and we postpone it to the remarks after the proof of proposition~\ref{prop:basic}. Thus altogether theorem~\ref{thm:main} gives an essentially complete picture except for the case (e$_1$) of finite quasisimple groups, which would require a close analysis of the representations of finite groups of Lie type generalizing the methods of Guralnick, Magaard, Malle and Tiep.

\medskip

For the sake of brevity, in what follows the term {\em group} will always be taken to include super groups. However, until section~\ref{sec:Lie} the term {\em dimension} will still refer to the ordinary dimension (as opposed to the super dimension).

\medskip

\section{Clifford-Mackey theory}

Let us say that $V\in \Rep_k(G)$ is {\em strongly irreducible} if for any non-central normal subgroup~$H\unlhd G$ of finite index the restriction $V|_H$ is irreducible.

\begin{prop} \label{irred_or_torus}
For any faithful $\epsilon$-small representation $V \in \Rep_k(G)$ one of the following cases occurs: \smallskip
\begin{enumerate}
 \item The representation $V$ is strongly irreducible. \smallskip
 \item $V$ is a $2$-homogenous monomial representation. In this case $\epsilon = -1$ or $V$ has dimension $\dim_k(V) \leq 2$. 
 \smallskip
 \item There exists an embedding $G \hookrightarrow \GO_4(k)$ such that $V$ is the restriction of the four-dimensional orthogonal standard representation.
\end{enumerate}
\end{prop}

{\em Proof.} Let $H\unlhd G$ be a normal subgroup. If the restriction~$V|_H$ is not isotypic, let $V|_H = W_1 \oplus \cdots \oplus W_n$ be its isotypic decomposition. Then $V\cong \Ind_{H_1}^G(W_1)$ is induced from a representation of the stabilizer~$H_1\leq G$ of $W_1$, and we get a splitting in two $G$-stable summands \bigskip
\[
 T_\epsilon(V)  \; \cong \; \Ind_{H_1}^G(T_\epsilon(W_1)) \oplus \Bigl[ \bigoplus_{i \neq j} W_i \otimes W_j \Bigr ]_\epsilon
\]
where the subscript $\epsilon$ in the second summand indicates the $\epsilon$-eigenspace of the commutativity constraint which flips the two factors of the tensor product. Since in the non-isotypic case we have $n>1$, $\epsilon$-smallness implies that $\dim_k(W_1)=1$, and~$\epsilon=-1$ or $\dim_k(V)=n=2$. All $W_i$ have dimension one, so $V|_{H}$ splits as a sum of pairwise distinct characters. Now $G$ acts by conjugation on the set $X$ of these characters, and the kernel $K$ of this permutation representation of $G$ is a normal subgroup which is abelian since $V$ is faithful. So (b) holds.

\medskip

Now suppose that $V|_H$ is isotypic. Then as in~\cite[th.~25.9]{Dornhoff} there are projective representations $U_1, U_2$ of $G$ such that $V \cong U_1 \otimes U_2$, where the restriction $U_1|_H$ is irreducible and where every $h\in H$ acts as the identity on~$U_2$. Then
\[
 T_\pm (V) \; \cong \; (T_+(U_1)\otimes T_{\pm}(U_2)) \oplus (T_-(U_1) \otimes T_{\mp}(U_2)),
\]
and since $V$ is small, one of the summands $T_{\epsilon_1}(U_1)\otimes T_{\epsilon_2}(U_2)$ must have dimension at most one. By direct inspection this can happen only if either $d_i = \dim_k(U_i) = 1$ for some $i\in \{1,2\}$, or $d_1=d_2=2$. Now
$
 V|_H \cong U_1 \oplus \cdots \oplus U_1 = d_2\cdot U_1
$
so that for $d_1=1$ the group $H$ is contained in the center $Z(G)$, acting on $V$ via scalar matrices. For~$d_2=1$ the restriction $V|_H$ remains irreducible. For $d_1=d_2=2$ case (c) occurs since $U_1, U_2\in \Rep_k(H)$ extend to projective representations of the whole group $G$ whose image then is contained in the product of its center with the special orthogonal similitude group $\Gl_2(k) \otimes \Gl_2(k) \cong \GSO_4(k)$. \qed

\section{Reduction to the quasisimple case}

Next we study the strongly irreducible $V\in \Rep_k(G)$. Let us temporarily call a group {\em basic} if it is either quasisimple or a finite $p$-group for some prime~$p$. For a given group $G$
we consider the following normal subgroup $H\unlhd G$:

\smallskip

\begin{itemize}
\item If $G^0 \subseteq Z(G)$, then $G = Z(G) \cdot S$ for some finite normal subgroup $S \unlhd G$, and fixing such a subgroup let $H=F^*(S)$ be its generalized Fitting subgroup as defined in~\cite[sect.~31]{Aschbacher}. \medskip

\item Otherwise we take $H=[G^0, G^0]$ to be the derived group of the connected component. The theory of reductive groups then implies $G^0 = Z(G^0)\cdot H$. \smallskip
\end{itemize}
In both cases we can find a central isogeny $\tilde{H} = H_1 \times \cdots \times H_n \twoheadrightarrow H$ such that the image of each $H_i$ is normal in $G$. Choosing the labelling in a suitable way, we may furthermore assume that for each $i$ we have a central isogeny $\tilde{H}_i = (G_i)^{s_i} \twoheadrightarrow H_i$ for $s_i$ copies of a suitable basic group~$G_i$ and that the images of these $s_i$ copies are permuted transitively by the adjoint action of $G$.

\begin{prop} \label{prop:basic}
For any faithful $\epsilon$-small strongly irreducible $V\in \Rep_k(G)$ with dimension $\dim_k(V)>1$ one of the following cases occurs: \smallskip
\begin{enumerate}
\item The group $H$ is quasisimple.
\smallskip
\item $(H, V|_H) \cong (G_1\otimes G_1, W\boxtimes W)$ for some very small $W\in \Rep_k(G_1)$, $H$ flips the two factors, and we have $\epsilon = -1$. \smallskip
\item $H$ contains an extraspecial $G$-stable subgroup $E$ of order $p^{2n+1}$ for some prime~$p$ such that $V|_E$ is irreducible of dimension $p^n$. \smallskip
\item We have an embedding $G\hookrightarrow \GO_4(k)$ such that $V$ is the restriction of the four-dimensional standard representation.
\end{enumerate}
\end{prop}

{\em Proof.} We first claim that $H\not \subseteq Z(G)$. Indeed, for the finite group case recall that the generalized Fitting subgroup contains its own centralizer~\cite{Aschbacher}, so $H \subseteq Z(G)$ would imply $S=H$ and then $G=Z(G)$ would be abelian. In the infinite case where $G^0$ is not central, the strong irreducibility implies that $V|_{G^0}$ is irreducible so that the connected reductive group $G^0$ cannot be a torus. Thus indeed $H\not \subseteq Z(G)$.

\medskip

Hence we can  assume that the image of each $H_i$ in $G$  is a {\em non-central} subgroup by discarding any occuring central components and saturating the other components with the centre. Since
$V|_{\tilde{H}} \cong U_1 \boxtimes \cdots \boxtimes U_n$ with irreducible $U_i \in \Rep_k(H_i)$, we get~$n=1$ by strong irreducibility. Hence $\tilde{H}\cong (G_1)^s$ for $s=s_1$ and again we get a decomposition
$V|_{\tilde{H}} \cong W_1 \boxtimes \cdots \boxtimes W_s$
with irreducible  $W_i \in \Rep_k(G_1)$,
but now the adjoint action of $G$ permutes the $s$ factors $G_1$ transitively so that all~$W_i$ are isomorphic to a single $W\in \Rep_k(G_1)$. In the decomposition
\[
 T_\epsilon(V)|_H \;\cong\; \bigoplus_{r=0}^s \; T_{r,\epsilon} \quad \textnormal{with} \quad
 T_{r,\epsilon} \;\; = \!\!\! \bigoplus_{\epsilon_1 \cdots \epsilon_s = \epsilon \atop \# \{ i \mid \epsilon_i = + 1\} = r} T_{\epsilon_1}(W) \boxtimes \cdots \boxtimes T_{\epsilon_s}(W)
\]
each summand $T_{r,\epsilon}$ is stable under the action of $G$. By smallness it then follows that $s\leq 2$, and for~$s=2$ the conclusions of (b) or (d) hold.

\medskip

So we may assume $s=1$ and $H=G_1$ is a basic group. If case (a) does not occur, then $H$ is a finite $p$-group for some prime $p$. Consider then a minimal $G$-stable non-central subgroup $M\unlhd H$. By minimality the subgroup $[M,M]$ is contained in~$A := M\cap Z(H)$ so that the quotient $U:=M/A$ is abelian. Looking at the \mbox{$p$-torsion} part of this quotient one obtains, again by minimality, that $U$ is elementary abelian. The commutator induces a bilinear map $[\cdot, \cdot]: U\times U \rightarrow A$, and if we identify~$A$ with a subgroup of $\bbG_m$ via Schur's lemma, $p\cdot U = 0$ implies that $[M, M]$ is contained in the subgroup $\mu_p \subseteq A$ of $p^\mathrm{th}$ roots of unity. So~$M/\mu_p$ is abelian and in fact elementary abelian: Otherwise by minimality  its $p$-torsion subgroup would lie in the cyclic group $A/\mu_p$ so that the abelian $p$-group $M/\mu_p$ would be cyclic as well. But then $M$ would be abelian, and this is impossible since it admits the faithful irreducible representation $V|_M$ of dimension $d>1$.

\medskip

Thus $M/\mu_p$ is elementary abelian, and we claim that the extraspecial case (c) occurs. Indeed, either $A=\mu_p$ or $A=\mu_{p^2}$. For $A=\mu_p$ the subgroup $E=M$ satisfies our requirements, so suppose that $A=\mu_{p^2}$. Since $M/\mu_p$ is elementary abelian, the Frattini subgroup is $\Phi(M)=\mu_p$ by~\cite[23.2]{Aschbacher}. The Frattini subgroup is the intersection of all maximal subgroups, so it follows that there exists a maximal subgroup $E \leq M$ which contains $\mu_p$ but not~$\mu_{p^2}$. Then $M=\mu_{p^2} \cdot E$, and $E\leq M$ is an extraspecial subgroup.
We will be done if we can show this subgroup is stable under the group $\Aut_A(M)$ of automorphisms of $M$ that are trivial on $A$. But this follows from the observation that every automorphism of $E$ which is trivial on $\mu_p$ extends uniquely to an element of $\Aut_A(M)$, which gives a natural identification $\Aut_A(M) \cong \Aut_{\mu_p}(E)$ compatible with the actions on $M$ and $E$. \qed

\medskip

We remark that the only instance of case (b) in proposition~\ref{prop:basic} with $\dim(H)>0$ is $G_1 \cong \Sl_m(k)$, acting on $W\cong k^m$ either via the standard representation or via its dual. Indeed this will follow from theorem~\ref{thm:main_thm} below, applied to the very small representation $W$ of the Lie algebra of $G_1$. Alternatively one could use~\cite{GT_Decompositions}. \medskip

In case (c) where $H$ contains a $G$-stable extraspecial $p$-group $E$, put $|E|=p^{1+2n}$ with $n\in \bbN$. For any non-trivial character $\omega: Z(E) \cong \bbZ/p\bbZ \to \bbG_m$ there exists a unique irreducible representation $V_\omega\in \Rep_k(E)$ of dimension $p^n$ on which~$Z(E)$ acts via the character $\omega$, and these are already all the irreducible representations of dimension $>1$ by~\cite[th.~31.5]{Dornhoff}. Hence in case (c) we have $V|_E \cong V_\omega$ for a uniquely determined character $\omega$. To decide which of the occuring representations are small, note that for the finite group $S$ such that $H=F^*(S)$, we have a natural homomorphism $S\to \Out(E)$. We now distinguish two cases depending on $p$.

\medskip

For $p>2$ we have $\omega^2\neq \one$, so $T_\epsilon(V)|_E$ is an isotypic multiple of~$V_{\omega^2}$. Then Mackey theory~\cite[th.~25.9]{Dornhoff} gives a tensor product decomposition
$T_\epsilon(V) \cong U \otimes W_\epsilon$
where~$U$ and~$W_\epsilon$ are projective representations of the group~$S$ such that $U|_E \cong V_{\omega^2}$ and such that every element of $E$ acts trivially on~$W_\epsilon$. Via the nondegenerate alternating bilinear form defined by the commutator on $E/Z(E) \cong (\bbF_p)^{2n}$ we can identify the image of $S$ in $\Out(E)$ with a subgroup of the symplectic group $\Sp_{2n}(\bbF_p)$, and one checks that in these terms $W_\epsilon$ is the restriction of the Weil representation~\cite{Gerardin_Weil} with dimension $(p^n + \epsilon)/2$. Accordingly $V$ is $\epsilon$-small iff the image of $S$ inside $\Sp_{2n}(\bbF_p)$ acts irreducibly on this Weil representation. 

\medskip

For $p=2$ on the other hand $\omega^2 = \one$, so that the restriction $T_\epsilon(V)|_E$ is a sum of characters. By~\cite{Winter_Automorphism} we can identify $\Out(E)$ with an orthogonal group $O^\pm_{2n}(\bbF_2)$ where the type $\pm$ of the quadratic form depends on $E$. Recall that a nondegenerate quadratic form on $(\bbF_2)^{2n}$ has type $\pm$ iff there are precisely $2^{n-1}(2^n \pm 1)$ isotropic vectors for this form. One then obtains the following identifications:

\begin{itemize}
\item If the quadratic form has $+$ type, the isotropic vectors in $(\bbF_2)^{2n}$ correspond precisely to the characters in $T_+(V)|_E$. \smallskip
\item If the quadratic form has $-$ type, the isotropic vectors in $(\bbF_2)^{2n}$ correspond precisely to the characters in $T_-(V)|_E$.
\end{itemize}
A similar interpretation holds for the anisotropic vectors. Hence it follows that~$V$ is small iff the image of $S$ inside $O^\pm_{2n}(\bbF_2)$ acts transitively on the nonzero isotropic resp.~anisotropic vectors. Note that the set of isotropic vectors always includes the zero vector as a single orbit, corresponding to the trivial summand $\one \hookrightarrow T_\epsilon(V)$.

\section{Lie super algebras} \label{sec:Lie}

It remains to determine all small $V\in \Rep_k(G)$ when $H=[G^0, G^0]$ is quasisimple and $V|_H$ is irreducible. By the classification of reductive super groups in~\cite{Weissauer_SS}, the Lie super algebra $\frakg$ of $H$ must then either be an ordinary simple Lie algebra or an orthosymplectic Lie super algebra $\osp_{1|2m}(k)$ with $m\in \bbN$. Note that $\Rep_k(H)$ is a full subcategory of $\Rep_k(\frakg)$, where the latter denotes the category of all Lie algebra representations of the Lie super algebra $\frakg$ on finite-dimensional super vector spaces over $k$. In particular~$V|_H$ defines an irreducible representation of $\frakg$. 

\medskip

The passage to representations of Lie algebras leads to a seemingly weaker notion of smallness. By the comments after theorem~\ref{thm:main} we know that $G/(G^0 \cdot Z_G(G^0))$ is a subgroup of $\Out(G^0)$ such that conjugation by any element $\varphi$ of this subgroup fixes the isomorphism type of $V|_H$. For an irreducible summand $W\hookrightarrow T_\epsilon(V)$ in~$\Rep_k(G)$ it may happen that the restriction $W|_H$ splits into several irreducible summands, but all these summands must be conjugate via automorphisms $\varphi$ as above. Abstracting from this situation, let us now denote by $\frakg$ any ordinary simple Lie algebra or $\osp_{1|2m}(k)$ with $m\in \bbN$. We say that a representation~$V\in \Rep_k(\frakg)$ is~{\em $\epsilon$-small} if either
$T_\epsilon(V) \cong W$ or $T_\epsilon(V) \cong W\oplus \one$,
where $W$ is a sum of irreducible representations which are all conjugate to each other via automorphisms $\varphi \in \mathrm{Aut}(\frakg)$ fixing the isomorphism type of $V$. To finish the proof of theorem~\ref{thm:main} we classify all irreducible small representations in this sense. For a uniform treatment the terms dimension, vector space, trace and Lie algebra will from now on be taken in the super sense for~$\osp_{1|2m}(k)$ but in the ordinary sense otherwise. 

\medskip

We denote by $\varpi_1, \dots, \varpi_m$ the fundamental dominant weights of $\frakg$ with respect to some fixed system of simple positive roots; see~\cite[sect.~2.1]{RS_Remarkable} for the orthosymplectic Lie algebra $\frakg = \osp_{1|2m}(k)$ whose Dynkin type we abbreviate by $BC_m$. Put
\[
 \beta_i \;=\; 
 \begin{cases}
 2\varpi_m & \textnormal{if $\frakg = \osp_{1|2m}(k)$ and $i=m$}, \\
 \;\, \varpi_i & \textnormal{otherwise}.
 \end{cases}
\]
The irreducible finite-dimensional representations of $\frakg$ are parametrized by highest weights $\lambda = \sum_{i=1}^m a_i \beta_i$ with $a_i\in \bbN_0$, see~\cite[th.~6]{Dj2}. For any such $\lambda$ we denote by~$V_\lambda$ the associated positive-dimensional irreducible representation. Note that in the super case negative-dimensional irreducible representations are obtained by the parity flip $W_\lambda = \Pi V_\lambda$ with $\dim(W_\lambda)=-\dim(V_\lambda)$ and $S^2(W_\lambda)\cong \Lambda^2(V_\lambda)$.

\setlength{\dashlinedash}{0.1pt}
\setlength{\dashlinegap}{1.5pt}

\begin{table}
\[
\begin{array}{|l:r|r|c|c|} \hline
 && \lambda & \epsilon = +1 & \epsilon = -1 \\ \hline
A_m & m\geq 1 & \beta_1, \, \beta_m & \star & \star \\ \cdashline{2-5}
& m=1 & 2\beta_1 & \circ & \star \\ 
 &   & 3\beta_1 & - & \circ \\ \cdashline{2-5}
 & m\geq 2 & 2\beta_1, \, 2\beta_m & - & \star \\ \cdashline{2-5}
 & m=3 & \beta_2 & \circ & \star \\ \cdashline{2-5}
 & m\geq 4 & \beta_2, \, \beta_{m-1} & - & \star \\ \cdashline{2-5}
 & m=5 & \beta_3 & - & \circ \\ \hline
B_m & m\geq 2 & \beta_1 & \circ & \star \\ \cdashline{2-5}
 & m=2 & \beta_2 & \star & \circ \\ \cdashline{2-5}
 & m=3 & \beta_3 & \circ & - \\ \hline
C_m & m\geq 3 & \beta_1 & \star & \circ \\ \cdashline{2-5}
 & m=3 & \beta_3 & - & \circ \\ \hline
D_m & m\geq 4 & \beta_1 & \circ & \star \\ \cdashline{2-5}
 & m=4 & \beta_3, \, \beta_4 & \circ & \star \\ \cdashline{2-5}
 & m=5 & \beta_4, \, \beta_5 & - & \star \\ \cdashline{2-5}
 & m=6 & \beta_5, \, \beta_6 & - & \circ \\ \hline
BC_m & m\geq 1 & \beta_1 & \star & \circ \\ \hline
E_6  && \beta_1, \, \beta_6 & - & \star \\ \hline
E_7 && \beta_7 & - & \circ \\ \hline
G_2 && \beta_1 & \circ & - \\ \hline
\end{array}
\]
\caption{The highest weights $\lambda$ of the small representations. The label $\star$ means $T_\epsilon(V_\lambda)$ is irreducible, $\circ$ means $T_\epsilon(V_\lambda) = W\oplus \one$ with $W$ irreducible, and the label $-$ means that $V_\lambda$ is not $\epsilon$-small.
}
\label{tab:smallreps}
\end{table}

\begin{thm} \label{thm:main_thm}
A positive-dimensional irreducible representation $V_\lambda \in \Rep_k(\frakg)$ is $\epsilon$-small iff its highest weight $\lambda$ appears in table~\ref{tab:smallreps}. 
\end{thm}

Note that for $\frakg = \fraksl_2(k)$ with its two-dimensional standard representation $st$, any irreducible representation is a symmetric power~$V_\lambda=S^n(st)$ of weight~$\lambda = n \beta_1$ for some $n\in \bbN$. In this case theorem~\ref{thm:main_thm} holds by direct inspection. A similar argument also works for $\frakg = \osp_{1|2}(k)$. Here we know from~\cite[th.~7 and th.~11]{Dj2} that for $\lambda = n\beta_1$ the even subalgebra~$\frakg_0 = \fraksl_2(k) \subset \frakg$ acts on~$V_\lambda=V_0 \oplus V_1$ via~$V_0 = S^n(st)$ and $V_1 = S^{n-1}(st)$. A short computation yields the action on the even and odd parts of the tensor square $T_\epsilon(V)$ and theorem~\ref{thm:main_thm} also holds in this case. Note that $\dim(V)=1$ for all irreducible representations $V$ of $\osp_{1|2}(k)$. For all other cases we have the following

\begin{lem} \label{lem:small_reps} 
For $\frakg \neq \osp_{1|2}(k)$ one has $\dim(V_\lambda)\leq \dim(\frakg)$ iff the highest weight $\lambda$ appears among those listed in tables~\ref{tab:replist1} or \ref{tab:replist2}.
\end{lem}

{\em Proof.} See~\cite{AEV} for the ordinary case. For~$\frakg = \osp_{1|2m}(k)$ with $m\geq 2$ we use the Kac-Weyl formula in~\cite[eq.~11]{CT_Supercharacters}. We embed the root system $BC_m$ into a Euclidean space with standard basis $\epsilon_1, \dots, \epsilon_m$ such that $\beta_i = \epsilon_1 + \cdots + \epsilon_i$ for all $i$. The irreducible representations of $\osp_{1|2m}(k)$ are parametrized by weights which in our basis are written $\lambda = (\lambda_1, \dots, \lambda_m)$ with integers $\lambda_1\geq \cdots \geq \lambda_m \geq 0$. The Kac-Weyl formula gives 
\[
 \dim(V_\lambda) \; \; \;=\;
 \prod_{1\leq i<j\leq m} \Bigl( \frac{\lambda_i - \lambda_j}{j-i}+1\Bigr)
 \; \; \cdot
 \prod_{1\leq i<j\leq m}
 \Bigl( \frac{\lambda_i + \lambda_j}{2m+1-i-j} + 1 \Bigr).
\]
For $\lambda_1\geq 2$ the second product is $\geq 2$. Then the classical Weyl formula for the first product shows that $\dim(V_\lambda)$ is at least twice the dimension of the irreducible representation of $\fraksl_m(k)$ with highest weight $\mu = (\lambda_1-\lambda_m, \dots, \lambda_{m-1}-\lambda_m)$. Using that $\dim(\fraksl_m(k)) \geq 2\dim(\osp_{1|2m}(k))$, it follows that $\mu$ is in the list for $A_{m-1}$ in table~\ref{tab:replist1}.
Since $\lambda = \mu + \lambda_m \cdot \beta_m$ and since increasing the weight by~$\beta_m$ increases the dimension, this leaves only finitely many cases. For $\lambda_1 = 1$ we have $\lambda = \beta_r$ with~$r\leq m$, and $\dim(V_\lambda) = {2m\choose r} - {2m\choose r-1}$ by the description in~\cite[sect.~5]{Dj2}.
\qed

\begin{cor} \label{cor:dimension_2} 
For $\frakg \neq \fraksl_2(k), \osp_{1|2}(k)$ and all weights $\lambda$ one has $\dim(V_\lambda)\geq 2$, with equality holding only in the single case~$(\frakg, \lambda)=(\osp_{1|4}(k), \beta_2)$.
\end{cor}

\section{Proof of theorem~\ref{thm:main_thm}} \label{sec:index}

Recall that $\frakg$ admits a unique invariant nondegenerate bilinear form~$(\cdot, \cdot)$ up to multiplication by a scalar~\cite[p.~94]{Sc}. Fixing any such form, we associate to any root~$\alpha$ a coroot $\alpha^\vee = 2\alpha / (\alpha, \alpha)$. Let~$\alpha_1, \dots, \alpha_m$ be a system of simple positive roots so that the fundamental weights $\varpi_i$ satisfy $( \alpha_i^\vee, \varpi_j) = \delta_{ij}$. Then $\rho = \varpi_1+\cdots +\varpi_m$ is half the sum of all positive roots, with the sign convention of~\cite{CT_Supercharacters}. For the proof of theorem~\ref{thm:main_thm} we consider the {\em index} of a representation $\varphi: \frakg\longrightarrow \frakgl(V)$, i.e.~the scalar $l(V)$ defined by $\mathrm{tr}(\varphi(X)\circ \varphi(Y)) = l(V)\cdot ( X, Y)$.

\begin{lem} \label{lem:index}
The index has the following properties. 
\begin{enumerate}
\item For the symmetric or alternating square of a representation $V$ it is given by the formula $l(T_\epsilon(V)) = (\dim(V) + 2\epsilon)\cdot l(V)$.
\smallskip

\item There exists a constant $\kappa \neq 0$ such that
$
\kappa \cdot l(V_\mu) =  \dim(V_\mu) \cdot c(\mu)$ for the scalar $c(\mu) = (\mu, \mu) + 2(\mu, \rho ) > 0
$ and for any highest weight $\mu\neq 0$. \smallskip

\item The index satisfies $l(\one)=0$, and it is invariant under automorphisms and additive for direct sums in the sense that $l(V\oplus V')= l(V) + l(V')$.
\end{enumerate}
\end{lem}

{\em Proof.} For (a) note that upon applying any tensor construction to~$V$ the index is multiplied by a constant depending only on $n=\dim(V)$. To compute this constant for $T_\epsilon(V)$, recall from~\cite[p.~128]{Sc} that~$\mathfrak{sl}(V)$ is simple for $n\neq 0$. It then only remains to check that $\mathrm{tr}((T_\epsilon(X))^2) = (n+2\epsilon)\, \mathrm{tr}(X^2)$ for a suitably chosen elementary matrix $X \in \mathfrak{sl}(V)$. For (b) one checks, by looking at the action on a highest weight vector, that the Casimir operator acts on $V_\mu$ by some fixed multiple of $c(\mu)$. The setting for $\osp_{1|2m}(k)$ is described in~\cite[p.~28]{Dj2} \cite[p.~223]{Dj1}. One then has $\kappa = \dim(Ad)\cdot c(Ad)$ for the adjoint representation $Ad$. Part (c) is obvious. \qed

\medskip

Via these index computations, we may now complete the classification of $\epsilon$-small representations for $\frakg \neq \fraksl_2(k), \osp_{1|2}(k)$ as follows.

\medskip

{\em Proof of theorem~\ref{thm:main_thm}}. Suppose that $V_\lambda$ is $\epsilon$-small. By corollary \ref{cor:dimension_2} we may assume that $n = \dim(V_\lambda) > 2$. Put $T_\epsilon(V_\lambda) = W \oplus \one^\delta$ where $\delta \in \{0,1\}$ denotes the multiplicity with which the trivial representation enters. Note that by smallness all highest weights $\mu$ occuring in $W$ are conjugate to each other. For any such~$\mu$ lemma~\ref{lem:index}(b),(c) hence imply that $\kappa \cdot l(W) =\dim(W) \cdot c(\mu) = ( n(n+\epsilon)/2 - \delta ) \cdot c(\mu)$ and $\kappa \cdot l(V_\lambda) = n \cdot c(\lambda)$. So lemma~\ref{lem:index}(a) shows
\begin{equation} \label{eq:casimir} \tag{$\star$}
 (n+2\epsilon) \cdot n \cdot c(\lambda) \;=\; \frac{1}{2} \Bigl(n(n+\epsilon) - 2\delta \Bigr) \cdot c(\mu).
\end{equation}
Now we distinguish between the symmetric and the alternating square. For $\epsilon = +1$ we may take $\mu = 2\lambda$. Then $c(\mu) = 4|\lambda|^2 + 4(\lambda, \rho)$. Since $c(\lambda)=|\lambda|^2 + 2(\lambda, \rho)$, equation~\eqref{eq:casimir} easily gives
\[
	(n-2\delta)\cdot \vert \lambda \vert^2 \;=\; 2\,( \lambda, \rho)
	\quad \textnormal{and hence} \quad 
	|\lambda | \; \leq \; \frac{2 |\rho|}{n-2\delta}
\]
by the Cauchy-Schwartz inequality. Let $\Delta_0$ be the set of simple positive roots of the even subalgebra $\frakg_0$. Then 
\[
 | (\lambda, \alpha^\vee ) |
 \; \leq \;
  |\lambda|  \cdot |\alpha^\vee|  
 \; \leq \;
 \frac{ 2 \, |\rho | \, |\alpha^\vee |}{n-2\delta}  
 \; < \; 
 \frac{ \dim( \frakg ) -1 }{n-2\delta} 
 \quad \textnormal{for any} \quad \alpha \in \Delta_0,
\]
where for the last inequality we have used the numerical values of $|\rho|^2$ and $R$ in table~\ref{tab:numerical} and our assumption $\frakg \neq \fraksl_2(k), \osp_{1|2}(k)$. On the other hand $(\lambda, \alpha^\vee)\in \bbZ$ for all $\alpha \in \Delta_0$, and for $\lambda \neq 0$ at least one of these scalar products is nonzero. Thus we can find $\alpha\in \Delta_0$ with $|(\lambda, \alpha^\vee )|\geq 1$. This implies $n - 2\delta < \dim(\frakg ) - 1$. Hence $\lambda$ is one of the highest weights in tables~\ref{tab:replist1} and \ref{tab:replist2} by lemma \ref{lem:small_reps}.

\medskip

It remains to discuss the case $\epsilon = -1$. By smallness all highest weights in $\Lambda^2(V_\lambda)$ are conjugate to each other via automorphisms fixing $\lambda$. Hence remark~\ref{rem:weights} below implies 
\begin{equation} \label{eq:lambda} \tag{$\star \star$}
 \lambda \;=\; r\cdot (\beta_{i_1} + \cdots + \beta_{i_s}) \quad \textnormal{for some $r\in \bbN$ and $i_1 < i_2 < \cdots < i_s$},  
\end{equation}
and that for all $i\in \{i_1, \dots, i_s\}$ the weight $\mu = 2\lambda - \alpha_i$ occurs as a highest weight in $\Lambda^2(V_\lambda)$. In what follows we fix~$i\in \{i_1, \dots, i_s\}$ with the smallest norm $|\beta_i|$. Since the norm of any simple positive root is given by the formula $|\alpha_i|^2 = 2\, (\alpha_i, \rho)$, we have $c(\mu)/2=c(\lambda) + |\lambda|^2 - 2(\lambda, \alpha_i)$ so that~\eqref{eq:casimir} becomes
\[ (n + 2\epsilon ) \cdot n \cdot c(\lambda) = \bigl(n(n+\epsilon) - 2\delta \bigr) \cdot \bigl( c(\lambda) + |\lambda|^2 - 2(\lambda, \alpha_i) \bigr). \]
Now for~$\epsilon = -1$ the first of the two factors on the right is $> (n+2\epsilon) \cdot n$ since by assumption $n > 2$ and $\delta \in \{0, 1\}$. Hence 
\[
 c(\lambda) + |\lambda|^2 - 2(\lambda, \alpha_i) \;<\; c(\lambda)
\]
and therefore
$
 2\, (\lambda, \alpha_i) > |\lambda|^2 \geq r^2 \cdot |\beta_i|^2 \cdot s,
$
where the second inequality comes from~\eqref{eq:lambda} together with the fact that all scalar products between~$\beta_{i_1}, \dots, \beta_{i_s}$ are nonnegative and $\beta_i$ has the smallest norm among all these weights. On the other hand $2\, (\lambda, \alpha_i)  = r \cdot 2\, (\beta_i, \alpha_i) = r \cdot |\alpha_i|^2$ by~\eqref{eq:lambda}. Hence $r_i  :=  \vert \alpha_i\vert^2 / \vert \beta_i\vert^2 > r\cdot s \geq 1$, which leaves only finitely many cases in view of table~\ref{tab:ri}. Note that for the Dynkin type~$A_m$ we may by duality assume $i<(m+1)/2$ so that $r_i \leq 4/i$. 
\qed

\medskip

For the convenience of the reader we include a proof of the following basic fact used in the above argument; see also~\cite[th.~5]{Aslaksen}.

\begin{rem} \label{rem:weights}
Let $\lambda = \sum_{i=1}^m a_i \beta_i$ with $a_i \in \bbN_0$. If $a_i > 0$, then the weight $2\lambda - \alpha_i$ appears as a highest weight in the alternating tensor square $\Lambda^2(V_\lambda)$.
\end{rem}

{\em Proof.} Let $v$ be a highest weight vector of $V_\lambda$. For $a_i > 0$ let $X_\pm \in \frakg_{\pm \alpha_i}$ be generators for the root spaces of the roots $\pm \alpha_i$ of $\frakg$ and put $H=[X_+, X_-]$. It then follows from $X_+v = 0$ that $X_+X_- v = H  v = (\alpha_i, \lambda) \cdot v \neq 0$. Since $v$ and $X_-v$ have different weights $\lambda$ resp.~$\lambda-\alpha_i$, this implies that $v\wedge X_-v \in \Lambda^2(V_\lambda)$ is a nonzero highest weight vector of weight $2\lambda - \alpha_i$. \qed

\section{An application to Brill-Noether sheaves} \label{sec:BN}

In this independent section we briefly discuss an application of theorem~\ref{thm:main} to algebraic geometry. Let $A$ be a complex abelian variety and $\bfD=\Dbc(A, \bbC)$ the derived category of bounded constructible sheaf complexes on it. Then the group law $a: A\times A \to A$ defines a convolution product $K*L = Ra_*(K\boxtimes L)$ on $\bfD$ which has all the formal properties of the tensor product in a Tannakian category, except that~$\bfD$ is not abelian. Passing to a certain abelian quotient category of perverse sheaves  one obtains for any $K\in \bfD$ a Tannaka group $G(K)$, see~\cite{KW_Vanishing}. 

\medskip

Now consider the special case where $A=\mathrm{Jac}(C)$ is the Jacobian variety of a smooth complex projective curve $C$ of genus $g>1$. Let $i: C\hookrightarrow A$ be a translate of the Abel-Jacobi embedding, and consider the corresponding constant perverse sheaf $\delta_C = i_*(\bbC_C[1]) \in \bfD$ on the image. The Tannaka group $G=G(\delta_C)$ depends on the choice of the embedding $i$. So, in what follows we normalize the embedding such that the highest exterior convolution power $\Lambda^{*(2g-2)}(\delta_C)$ is represented by the skyscraper sheaf $\one$ of rank one supported in the origin. This is possible because by prop.~10.1 of loc.~cit.~this highest exterior convolution power is a skyscraper sheaf of rank one. Now, assuming $g>2$, the classification in theorem~\ref{thm:main} allows to give a rather short proof of the following result of~\cite{Weissauer_BN}.

\begin{thm}
Let $C$ be a smooth complex projective curve of genus $g>1$ which is embedded into its Jacobian variety $A=\mathit{Jac}(C)$ as above. Then
\[
 G(\delta_C) \;=\; 
 \begin{cases}
 \, \Sl_{2g-2}(\bbC) & \textnormal{\em if $C$ is not hyperelliptic}, \\
 \Sp_{2g-2}(\bbC) & \textnormal{\em if $C$ is hyperelliptic}.
 \end{cases}
\]
\end{thm}

{\em Proof for $g>2$.} For hyperelliptic curves $C$ the Abel-Jacobi map $f: C^2 \rightarrow A$ is generically finite of degree two over its image, but blows down the hyperelliptic linear series $g_2^1$ to a point $a\in A(\bbC)$. By our normalization of the embedding $C\hookrightarrow A$ we can assume $a=0$. One then checks that $\delta_C * \delta_C = Rf_*(\bbC_{C^2}[2]) = \delta_+ \oplus \delta_- \oplus \one$ for certain simple perverse sheaves $\delta_\pm$ and the rank one skyscraper sheaf $\one$ with support in the origin.
The definition of the symmetry constraint in~\cite{Weissauer_BN} shows that~$\one$ lies in the {\em alternating} convolution square of $\delta_C$. If $G=G(\delta_C)$ denotes our Tannaka group and if $V\in \Rep_k(G)$ denotes the representation corresponding to the perverse sheaf $\delta_C$, it follows that the symmetric square $T_+(V)$ is irreducible and that $T_-(V)$ decomposes into an irreducible plus a trivial representation.

\medskip

The $\epsilon$-smallness of $V$ for $\epsilon = +1$ rules out case (b) in theorem~\ref{thm:main}. Case (d) is ruled out for the same reason because by~\cite{KW_Vanishing} the dimension of any representation of~$G$ is the Euler characteristic of the underlying perverse sheaf, which in our situation is $d=\dim_k(V)=2g-2>2$ for $g>2$. Since~$T_+(V)$ is irreducible whereas the symmetric square of the standard representation of the orthogonal group is not, case (c) is impossible. 
Case~(e) cannot occur either since the group of connected components of the Tannaka group of any perverse sheaf is abelian~\cite{Weissauer_Connected}. So we must be in case (a), and we are looking for entries in table~\ref{tab:smallreps} with the label $\star$ for $\epsilon = +1$ and $\circ$ for $\epsilon = -1$. Since we are dealing with ordinary groups, the only possibility is the standard representation of $\Sp_{2m}(k)$ where $2m=d=2g-2$; for $g=3$ the isomorphism $B_2\cong C_2$ must be used in reading the table. The non-hyperelliptic case is similar but here no summand $\one$ occurs. \qed



\setlength{\dashlinedash}{0.1pt}
\setlength{\dashlinegap}{1.5pt}

\begin{table}
\footnotesize
\[
\begin{array}{|l:r|r|r|r|} \hline  \medrowheight
  &  & \lambda & S^2(V_\lambda) & \Lambda^2(V_\lambda) 
  \\ \hline 
  A_m & m = 1 & \beta_1 & V_{2\beta_1} & \one
	\\ \cdashline{2-5} 
	& m\geq 2 & \beta_1 & V_{2\beta_1} & V_{\beta_2}
	\\ \cdashline{3-5} 
	&& \beta_m & V_{2\beta_m} & V_{\beta_{m-1}}
	\\ \cdashline{2-5} 
	& m\geq 2	& 2\beta_1 &  V_{4\beta_1}\oplus V_{2\beta_2} & V_{2\beta_1+\beta_2}
	 \\ \cdashline{3-5} 
	&& 2\beta_m & V_{4\beta_m}\oplus V_{2\beta_{m-1}} & V_{2\beta_m + \beta_{m-1}}
	\\ \cdashline{2-5} 
	&  m \geq 4 & \beta_2 & V_{2\beta_2}\oplus V_{\beta_4} & V_{\beta_1 + \beta_3}
	\\ \cdashline{3-5} 
	&& \beta_{m-1} & V_{2\beta_{m-1}}\oplus V_{\beta_{m-3}} & V_{\beta_m+\beta_{m-2}}
	\\ \cdashline{2-5} 
	&  m=3 & \beta_2 & V_{2\beta_2}\oplus \one & V_{\beta_1+\beta_3}
	\\ \cdashline{2-5}
	&  m= 5 & \beta_3 & V_{2\beta_3}\oplus V_{\beta_1+\beta_5} & V_{\beta_2+\beta_4} \oplus \one
	\\ \cdashline{2-5} 
	& m= \;\;\; & \beta_3 & V_{2\beta_3}\oplus V_{\beta_1+\beta_5} & V_{\beta_2+\beta_4} \oplus V_{\beta_6}
	\\ \cdashline{3-5} 
	& 6,7 &  \beta_{m-2} & V_{2\beta_{m-2}}\oplus V_{\beta_m+\beta_{m-4}} & V_{\beta_{m-1}+\beta_{m-3}} \oplus V_{\beta_{m-5}}
	\\ \hline  
B_m & m\geq 2 & \beta_1 & V_{2\beta_1}\oplus \one & V_{\beta_2}
	\\ \cdashline{2-5}
	& m = 2 & \beta_2 & V_{2\beta_2} & V_{\beta_1} \oplus \one
	\\ \cdashline{2-5}
	& m=3 & \beta_3 & V_{2\beta_3}\oplus \one & V_{\beta_1}\oplus V_{\beta_2}
	\\ \cdashline{2-5} 
	& m=4 & \beta_4 & V_{2\beta_4}\oplus V_{\beta_1}\oplus \one & V_{\beta_2}\oplus V_{\beta_3}
	\\ \cdashline{2-5}
	& m=5 & \beta_5 & V_{2\beta_5}\oplus V_{\beta_2} 
	\oplus V_{\beta_1} & V_{\beta_3} \oplus V_{\beta_4} \oplus \one
	\\ \cdashline{2-5} 
	& m=6 & \beta_6 & V_{2\beta_6}\oplus V_{\beta_3} 
	\oplus V_{\beta_2} & V_{\beta_1} \oplus V_{\beta_4} \oplus V_{\beta_5} \oplus \one
	\\ \hline 
C_m & m\geq 3 & \beta_1 & V_{2\beta_1} & V_{\beta_2} \oplus \one
	\\ \cdashline{2-5} 
	& m \geq 4 & \beta_2 & V_{\beta_4}\oplus V_{2\beta_2}\oplus V_{\beta_2} \oplus \one &   V_{2\beta_1} \oplus V_{\beta_1+\beta_3}
	\\ \cdashline{2-5} 
	& m = 3 & \beta_2 & V_{2\beta_2}\oplus V_{\beta_2}\oplus \one & V_{2\beta_1} \oplus V_{\beta_1+\beta_3}
	\\ \cdashline{3-5} 
	& & \beta_3 & V_{2\beta_3} \oplus V_{2\beta_1} & V_{2\beta_2} \oplus \one
	\\ \hline 
D_m & m\geq 4 & \beta_1 & V_{2\beta_1} \oplus \one & V_{\beta_2}
	\\ \cdashline{2-5} 
	& m = 4 & \beta_3 & V_{2\beta_3} \oplus \one & V_{\beta_2}
	\\ \cdashline{3-5} 
	& & \beta_4 & V_{2\beta_4} \oplus \one & V_{\beta_2}
	\\ \cdashline{2-5}
	& m = 5 & \beta_4 & V_{2\beta_{4}}\oplus V_{\beta_{1}} & V_{\beta_3}
	\\ \cdashline{3-5}
	& & \beta_5 & V_{2\beta_{5}}\oplus V_{\beta_{1}} & V_{\beta_3}
	\\ \cdashline{2-5} 
	& m = 6 & \beta_5 & V_{2\beta_{5}}\oplus V_{\beta_{2}} & V_{\beta_4} \oplus \one
	\\ \cdashline{3-5}
	& & \beta_6 & V_{2\beta_{6}}\oplus V_{\beta_{2}} & V_{\beta_4} \oplus \one
	\\ \cdashline{2-5} 
	& m = 7 & \beta_6 & V_{2\beta_{6}}\oplus V_{\beta_{3}} & V_{\beta_1} \oplus V_{\beta_5}
	\\ \cdashline{3-5} 
	& & \beta_7 & V_{2\beta_{7}}\oplus V_{\beta_{3}} & V_{\beta_1} \oplus V_{\beta_5}
	\\ \hline 
BC_m & m\geq 2 & \beta_1 & V_{2\beta_1} & V_{\beta_2} \oplus \one \\ \cdashline{2-5} 
	& m\geq 4 & \beta_2 & V_{2\beta_2} \oplus V_{\beta_2} \oplus V_{\beta_4} \oplus \one
	& V_{2\beta_1} \oplus V_{\beta_1 + \beta_3}
	\\ \cdashline{2-5} 
	& m=2 & \beta_1 + \beta_2 & V_{2\beta_1+2\beta_2} \oplus 2V_{2\beta_1+\beta_2} \oplus 2V_{2\beta_1} 
	& V_{4\beta_1} \oplus V_{2\beta_1 + \beta_2} \oplus V_{3\beta_2} \oplus 2V_{2\beta_2}  \\ 
	&&& \oplus W_{3\beta_1} \oplus W_{\beta_1+2\beta_2}  
	& \oplus 2V_{\beta_2} \oplus W_{3\beta_1 + \beta_2} \oplus W_{\beta_1+2\beta_2} \\  
	&&& \oplus 2 W_{\beta_1+\beta_2} 
	& \oplus W_{\beta_1 + \beta_2} \oplus W_{\beta_1}\oplus \one \\ \cdashline{3-5} 
	&  & \beta_2 & V_{2\beta_2} \oplus V_{\beta_2} \oplus W_{\beta_1} 
	& V_{2\beta_1} \oplus W_{\beta_1+\beta_2} \oplus \one \\ \cdashline{2-5} 
	& m=3 & \beta_2 & V_{2\beta_2} \oplus V_{\beta_2} \oplus W_{\beta_3} 
	& V_{2\beta_1} \oplus V_{\beta_1 + \beta_3}\oplus \one \\ \cdashline{3-5} 
	& & \beta_3 & V_{2\beta_1} \oplus V_{\beta_1 + \beta_3} \oplus V_{2\beta_3}  
	& V_{2\beta_2} \oplus V_{\beta_2} \\ 
	&&& \oplus W_{\beta_1 + \beta_2}
	& \oplus W_{\beta_1} \oplus W_{\beta_2 + \beta_3} \oplus W_{\beta_3}  \oplus \one\\ \cdashline{2-5} 
	& m=4 & \beta_4 & V_{2\beta_1} \oplus V_{\beta_2 + \beta_4} \oplus V_{\beta_2} \oplus V_{2\beta_4} 
	& V_{2\beta_2} \oplus V_{\beta_1+\beta_3} \oplus V_{2\beta_3} \\ 
	&&&  \oplus V_{\beta_4} \oplus W_{\beta_1} \oplus W_{\beta_2+\beta_3} 
	& \oplus W_{\beta_1+\beta_2} \oplus W_{\beta_1 + \beta_4}  \\ 
	&&& \oplus W_{\beta_3}\oplus \one 
	& \oplus W_{\beta_3 + \beta_4} 
	\\ \cdashline{2-5} 
	& m=5 & \beta_5 & V_{2\beta_1} \oplus V_{\beta_1+\beta_3} \oplus V_{\beta_1 + \beta_5} 
	& V_{2\beta_2} \oplus V_{\beta_2+\beta_4}  \oplus V_{\beta_2} \oplus V_{2\beta_4}  \\ 
	&&& \oplus V_{2\beta_3} \oplus V_{\beta_3+\beta_5} \oplus V_{2\beta_5} 
	& \oplus V_{\beta_4} \oplus W_{\beta_1} \oplus W_{\beta_2 + \beta_3}  \\ 
	&&& \oplus W_{\beta_1 + \beta_2} \oplus W_{\beta_1 + \beta_4} 
	& \oplus W_{\beta_2 + \beta_5}  \oplus W_{\beta_3}   \\ 
	&&& \oplus W_{\beta_3 + \beta_4}
	& \oplus W_{\beta_4+\beta_5} \oplus W_{\beta_5} \oplus \one \\ \hline	
E_6 &  & \beta_1 & V_{2\beta_1}\oplus V_{\beta_6} & V_{\beta_3}
	\\ \cdashline{3-5}
	&& \beta_6 & V_{2\beta_6}\oplus V_{\beta_1} & V_{\beta_5}
	\\ \hline 
E_7    & & \beta_7 & V_{2\beta_7}\oplus V_{\beta_1} & V_{\beta_6} \oplus \one
	\\ \hline 
F_4 && \beta_4 & V_{2\beta_4}\oplus V_{\beta_4} \oplus \one & V_{\beta_3} \oplus V_{\beta_1}
	\\ \hline 
G_2 && \beta_1 & V_{2\beta_1}\oplus \one & V_{\beta_1} \oplus V_{\beta_2} \\ \hline
\end{array}
\]
\caption{All $\lambda$ with $1 < \dim(V_\lambda) < \dim(\frakg)$. For $\frakg=\osp_{1|2m}(k)$ we denote the parity shifts of the highest weight modules by $W_\mu = \Pi V_\mu$. }
\label{tab:replist1}
\end{table}

\begin{table}
\footnotesize
\[
\begin{array}{|l:r|r|r|r|} \hline  \medrowheight
  && \lambda & S^2(V_\lambda) & \Lambda^2(V_\lambda) 
  \\ \hline 
A_m & m = 1 & 2\beta_1 & V_{4\beta_1} \oplus \one & V_{2\beta_1} \\ \cdashline{2-5} 
	& m = 2 & \beta_1 + \beta_2 
	& V_{2\beta_1 + 2\beta_2} \oplus V_{\beta_1 + \beta_2} \oplus \one
	& V_{3\beta_1} \oplus V_{3\beta_2} \oplus V_{\beta_1 + \beta_2} \\ 
	\cdashline{2-5} 
	& m\geq 2 & \beta_1 + \beta_m 
	& V_{2\beta_1 + 2\beta_m} \oplus V_{\beta_2 + \beta_{m-1}} 
	& V_{\beta_2 + 2\beta_m} \oplus V_{2\beta_1 + \beta_{m-1}}  \\ 
	&&& \oplus V_{\beta_1+\beta_m} \oplus \one & \oplus V_{\beta_1 + \beta_m}
	\\ \hline 
B_m 
	& m = 2 & 2\beta_2 
	&  V_{\beta_1} \oplus V_{2\beta_1} \oplus V_{4\beta_2} \oplus \one 
	&  V_{\beta_1+2\beta_2} \oplus V_{2\beta_2} \\ \cdashline{2-5} 
	& m =3 & \beta_2 
	& V_{2\beta_1} \oplus V_{2\beta_2} \oplus V_{2\beta_3} \oplus \one & V_{\beta_1 + 2\beta_3} \oplus V_{\beta_2} \\ \cdashline{2-5}
	& m = 4 & \beta_2
	& V_{2\beta_1} \oplus V_{2\beta_2} \oplus V_{2\beta_4} \oplus \one & V_{\beta_1 + \beta_3} \oplus V_{\beta_2}
	\\ \cdashline{2-5} 
	& m \geq 5 & \beta_2
	& V_{2\beta_1} \oplus V_{2\beta_2} \oplus V_{\beta_4} \oplus \one & V_{\beta_1 + \beta_3} \oplus V_{\beta_2}
	\\ \hline
C_m 
	& m\geq 3 & 2\beta_1
	& V_{4\beta_1} \oplus V_{2\beta_2} \oplus V_{\beta_2} \oplus \one & V_{2\beta_1} \oplus V_{2\beta_1 + \beta_2}
	\\ \hline 
D_m 
	& m = 4 & \beta_2
	& V_{2\beta_1} \oplus V_{2\beta_2} \oplus V_{2\beta_3} \oplus V_{2\beta_4} \oplus \one 
	& V_{\beta_2} \oplus V_{\beta_1+\beta_3+\beta_4} \\ \cdashline{2-5} 
	& m = 5 & \beta_2
	& V_{2\beta_1} \oplus V_{2\beta_2} \oplus V_{\beta_4+\beta_5} \oplus \one & V_{\beta_2} \oplus V_{\beta_1 + \beta_3}
	\\ \cdashline{2-5} \medrowheight
	& m \geq 6 & \beta_2
	& V_{2\beta_1} \oplus V_{2\beta_2} \oplus V_{\beta_4} \oplus \one & V_{\beta_2} \oplus V_{\beta_1 + \beta_3}
	\\ \hline \medrowheight
BC_m & m \geq 2 & 2\beta_1 & V_{4\beta_1} \oplus V_{2\beta_2} \oplus V_{\beta_2} \oplus \one  & V_{2\beta_1+\beta_2} \oplus V_{2\beta_1} \\ \cdashline{2-5} 
	& m = 4 & \beta_3 & V_{2\beta_1} \oplus V_{\beta_1+\beta_3}  \oplus V_{2\beta_3}
 	& V_{2\beta_2} \oplus V_{\beta_2} \oplus V_{\beta_2+\beta_4} \\ 
	&& &  \oplus W_{\beta_1+\beta_4}
	&  \oplus V_{\beta_4} \oplus W_{\beta_3} \oplus \one  	\\ 
 	\hline 
E_6 && \beta_2 & V_{2\beta_2} \oplus V_{\beta_1+\beta_6} \oplus \one & V_{\beta_2} \oplus V_{\beta_4}
	\\ \hline 
E_7 && \beta_1 & V_{2\beta_1} \oplus V_{\beta_6} \oplus \one & V_{\beta_1} \oplus V_{\beta3}
	\\ \hline 
E_8 && \beta_8 & V_{\beta_1} \oplus V_{2\beta_8} \oplus \one & V_{\beta_7} \oplus V_{\beta_8} 
	\\ \hline 
F_4 && \beta_1 & V_{2\beta_1} \oplus V_{2\beta_4} \oplus \one & V_{\beta_1} \oplus V_{\beta_2} 
	\\ \hline 
G_2 && \beta_2 & V_{2\beta_1} \oplus V_{2\beta_2} \oplus \one & V_{3\beta_1} \oplus V_{\beta_2} \\ \hline
\end{array}
\]
\caption{All $\lambda$ with $1 < \dim(V_\lambda) = \dim(\frakg)$. For the ordinary simple Lie algebras precisely the adjoint representations occur~\cite{AEV}.}
\label{tab:replist2}
\end{table}

\begin{table} 
\footnotesize
\[
\begin{array}{|l|c|c|c|l|l|} \hline \bigrowheight
	& |\rho|^2 & R & \dim(\frakg) & r_i \;\; \textnormal{for} \;\; i \;=\; 1,\dots, m& |\Out(\frakg)|
	\\ \hline \bigrowheight
	A_m & \frac{m(m+1)(m+2)}{12} & \sqrt{2} & m(m+2) & \frac{2(m+1)}{i(m+1-i)} & 2
	\\ \cdashline{1-6} \bigrowheight 
	B_m & \frac{m(2m-1)(2m+1)}{12} & 2 & m(2m+1) & \frac{2}{i} (1+\delta_{im}) & 1
	\\ \cdashline{1-6} \bigrowheight
	C_m & \frac{m(m+1)(2m+1)}{6} & \sqrt{2} & m(2m+1) & \frac{2}{i} (1+\delta_{im}) & 1
	\\ \cdashline{1-6} \bigrowheight
	D_m & \frac{m(m-1)(2m-1)}{6} & \sqrt{2} & m(2m-1) & \frac{2}{i} \hspace{0.83em} \textnormal{if $i<m-1$} & 2  \;\; \textnormal{if $m\neq 4$} \\ \bigrowheight
	&&&& \frac{8}{m} \hspace{0.5em} \textnormal{if $i\in \{m-1, m\}$} & 6 \;\; \textnormal{if $m=4$}
	\\ \cdashline{1-6} \bigrowheight
	BC_m & \frac{m(2m-1)(2m+1)}{12} & \sqrt{2} & m(2m-1) & \frac{2}{i} (1+\delta_{im}) & 1
	\\ \cdashline{1-6} \bigrowheight 
	E_6 & 78 & \sqrt{2} & 78 & \frac{3}{2}, 1, \frac{3}{5}, \frac{1}{3}, \frac{3}{5}, \frac{3}{2} & 2
	\\ \cdashline{1-6} \bigrowheight
	E_7 & \frac{399}{2} & \sqrt{2} & 133 & 1, \frac{4}{7}, \frac{1}{3},  \frac{1}{6},  \frac{4}{15}, \frac{1}{2}, \frac{4}{3} & 1
	\\ \cdashline{1-6} \bigrowheight
	E_8 & 620 & \sqrt{2} & 248 & \frac{1}{2}, \frac{1}{4}, \frac{1}{7}, \frac{1}{15}, \frac{1}{10}, \frac{1}{6}, \frac{1}{3}, 1 & 1
	\\ \cdashline{1-6} \bigrowheight
	F_4 & 39 & 2 & 52 & 1, \frac{1}{3}, \frac{1}{3}, 1 & 1
	\\ \cdashline{1-6} \bigrowheight
	G_2 & 14 & \sqrt{2} & 14 & 1 & 1
	\\ \hline
\end{array}
\]
\caption{Some numerical values. Here we write $r_i = |\alpha_i|^2/|\beta_i|^2$ and put $R = \max_{\alpha \in \Delta_0} |\alpha^\vee|$ for the set~$\Delta_0$ of simple positive roots of $\frakg_0$.}
\label{tab:numerical} \label{tab:ri}
\end{table}

\bibliographystyle{amsplain}
\bibliography{SmallTensorGenerator}

\end{document}